\documentclass[12pt,a4paper]{article}
\usepackage[T1]{fontenc}
\usepackage{lmodern}
\usepackage[margin=2.7cm]{geometry}
\usepackage{amsmath,amssymb,amsthm}
\usepackage{tikz}
\usetikzlibrary{calc}
\usepackage[hidelinks]{hyperref}
\hypersetup{pdftitle={Unique-maximum colourings of planar graphs with respect to
neighbourhoods},pdfauthor={Jelena Sedlar, Riste Skrekovski}}
\definecolor{gA}{RGB}{0,90,170}
\definecolor{gB}{RGB}{0,130,60}
\definecolor{gC}{RGB}{200,100,0}
\newtheorem{theorem}{Theorem}
\newtheorem{lemma}[theorem]{Lemma}
\newtheorem{corollary}[theorem]{Corollary}
\newtheorem{conjecture}[theorem]{Conjecture}
\theoremstyle{definition}
\newtheorem{definition}[theorem]{Definition}
\newcommand{\pUMc}{\chi_{\mathrm{pUMc}}}
\newcommand{\pUMo}{\chi_{\mathrm{pUMo}}}
\newcommand{\iUMc}{\chi_{\mathrm{iUMc}}}
\newcommand{\iUMo}{\chi_{\mathrm{iUMo}}}
\begin{document}

\begin{center}
{\Large\bfseries Unique-maximum colourings of planar graphs\\[0.1cm]
with respect to neighbourhoods}\\[0.4cm]
Jelena Sedlar$^{1,3}$, Riste \v{S}krekovski$^{2,3}$\\[0.2cm]
{\small $^{1}$ \textit{University of Split, Faculty of Civil Engineering, Architecture and
Geodesy, Split, Croatia}}\\[0.05cm]
{\small $^{2}$ \textit{University of Ljubljana, FMF, Ljubljana, Slovenia}}\\[0.05cm]
{\small $^{3}$ \textit{University of Novo Mesto, FIS, Novo Mesto, Slovenia}}
\end{center}

\begin{abstract}
A colouring of the vertices of a graph by positive integers is a
unique-maximum colouring with respect to open (resp.\ closed) neighbourhoods if
the maximum colour appears exactly once in the open (resp.\ closed)
neighbourhood of every vertex. Fabrici, Lu\v{z}ar, Rindo\v{s}ov\'a and Sot\'ak
(\emph{Discrete Appl.\ Math.}\ \textbf{324} (2023) 80--92) conjectured that four
colours suffice for closed neighbourhoods in planar graphs. We prove this
conjecture by colouring a suitable independent dominating set with the three
largest colours. They further conjectured that six colours suffice if the
colouring is also proper. We disprove this conjecture by a planar graph which
needs seven colours, and we show that seven colours always suffice. For open
neighbourhoods, we show that eight colours suffice, even for proper colourings,
which improves the previous bound of ten. In the proper case, our graph again
needs seven colours. If the colouring need not be proper and the graph has no
vertices of degree $2$, then six colours suffice.
\end{abstract}


\section{Introduction}

Conflict-free colourings of hypergraphs were introduced by Even, Lotker, Ron and
Smorodinsky \cite{Even}, and conflict-free colourings of graphs with respect to
neighbourhoods were studied in \cite{Abel, PachTardos}. We consider a related
notion. A colouring of the vertices of a graph $G$ by positive integers is a
\emph{unique-maximum colouring with respect to open} (resp.\ \emph{closed})
\emph{neighbourhoods} if, for every vertex $v$, the maximum colour appears on
exactly one vertex of the open neighbourhood $N(v)$ (resp.\ of the closed
neighbourhood $N[v]$). We call such a colouring an \emph{iUMo-colouring} (resp.\
an \emph{iUMc-colouring}). If the colouring is also proper, we call it a
\emph{pUMo-colouring} (resp.\ a \emph{pUMc-colouring}). The minimum numbers of
colours in these colourings are denoted by $\iUMo(G)$, $\iUMc(G)$, $\pUMo(G)$
and $\pUMc(G)$, respectively. For the open variants, we only consider graphs
without isolated vertices.

For a class of graphs $\mathcal{C}$ and an invariant $\chi_{\ast}$, we write
$\chi_{\ast}(\mathcal{C})=\max\{\chi_{\ast}(G):G\in\mathcal{C}\}$. For the open
variants, the maximum is taken over the graphs of $\mathcal{C}$ without isolated
vertices. For all invariants and classes in this note, the maximum exists by the
upper bounds below. We denote by $\mathcal{P}$ the class of planar graphs. For a
vertex $v$ of a graph $G$, we write $N_{G}(v)$ and $N_{G}[v]$ for its open and
closed neighbourhoods and $d(v)$ for its degree, and we omit the subscript when
the graph is clear from the context. The \emph{join} $A\vee B$ of two
vertex-disjoint graphs $A$ and $B$ is obtained from their disjoint union by
adding all edges between $A$ and $B$.

The four invariants were studied by Fabrici, Lu\v{z}ar, Rindo\v{s}ov\'a and
Sot\'ak \cite{FLRS}. Among other results, they proved that
\[
4\leq\iUMc(\mathcal{P})\leq6,\quad5\leq\iUMo(\mathcal{P})\leq10,\quad
6\leq\pUMc(\mathcal{P})\leq8,\quad6\leq\pUMo(\mathcal{P})\leq10.
\]
Their upper bounds combine a dominating colour class of a proper colouring with
the facial unique-maximum colouring theorem of Wendland \cite{Wendland}, which
they apply to a suitable facial closure of the graph. They proposed the
following conjectures.

\begin{conjecture}
\label{Con_FLRS}It holds that
(i) $\iUMc(\mathcal{P})=4$, (ii) $\iUMo(\mathcal{P})=5$, (iii) $\pUMc(\mathcal{P})=6$ and
(iv) $\pUMo(\mathcal{P})\leq7$.
\end{conjecture}

In this note we show that part (i) of Conjecture \ref{Con_FLRS} is true and
that part (iii) is false. Namely, we prove that
\[
\iUMc(\mathcal{P})=4\quad\text{and}\quad\pUMc(\mathcal{P})=7.
\]
The lower bound $\pUMc(\mathcal{P})\geq7$ is attained by a planar graph on $151$
vertices, which is obtained from $K_{4}$ by attaching seven copies of an
outerplanar graph from \cite{FLRS} (Theorem \ref{Tm_lower7}). Both upper bounds
follow from the iterated distance-three selector of Abel et al.\ \cite{Abel}, in
the form given recently by Jim\'enez, Lintzmayer and Sambinelli \cite{JLS}. The
key property is that in every closed neighbourhood the \emph{maximum} selector
colour is unique (Lemma \ref{Lemma_selector}). This is implicit in \cite{JLS},
where the newly chosen set always receives the largest colour. We also prove
that
\[
7\leq\pUMo(\mathcal{P})\leq8\quad\text{and}\quad\iUMo(\mathcal{P})\leq8,
\]
which improves the bounds of \cite{FLRS}. For planar graphs without vertices of
degree $2$, we obtain $\iUMo(G)\leq6$ (Theorem \ref{Tm_iUMo6}). For this, we
prove a version of a lemma of Fabrici and G\"{o}ring for plane multigraphs
(Lemma \ref{Lemma_FG}). Parts (ii) and (iv) of Conjecture \ref{Con_FLRS} remain
open. In the concluding section we explain why our method does not give $7$ for
$\pUMo$.

\section{Unique-maximum selectors}

We use the following unique-maximum version of the selectors of \cite{JLS}.

\begin{definition}
Let $G$ be a graph and let $p$ be a positive integer. A \emph{UM-selector of
width} $p$ in $G$ is a pair $(S,\varphi)$, where $S\subseteq V(G)$ is an
independent set and $\varphi:S\rightarrow\{1,\ldots,p\}$, such that for every
vertex $x\in V(G)$ the set $N[x]\cap S$ is non-empty and the maximum of
$\varphi$ on $N[x]\cap S$ is attained by exactly one vertex.
\end{definition}

UM-selectors of the components of a graph together form a UM-selector of the
whole graph. Our construction of UM-selectors is based on the following lemma
of Jim\'enez, Lintzmayer and Sambinelli, which isolates the distance-three
construction of Abel et al.\ \cite{Abel}. All our upper bounds depend on it, so
we include its short proof from \cite{JLS}.

\begin{lemma}
\label{Lemma_JLS}Every nonempty connected graph $G$ contains a set $D\subseteq
V(G)$ such that \emph{(L1)} distinct vertices of $D$ are at distance at least
$3$, and \emph{(L2)} if a graph $J$ is a minor of $R=G-N[D]$, then $K_{1}\vee J$
is a minor of $G$.
\end{lemma}

\begin{proof}
Start with $D=\{v\}$ for an arbitrary vertex $v$. As long as some vertex is at
distance at least $3$ from $D$, add to $D$ a vertex at distance exactly $3$ from
$D$; such a vertex lies on a shortest path from $D$ to a vertex at distance at
least $3$. Then (L1) holds. At the end, every vertex is at distance at most $2$
from $D$, so every vertex of $R$ has a neighbour in $N[D]$. When a vertex $w$ is
added, there is a path $waw^{\prime\prime}w^{\prime}$ with $w^{\prime}$ already in
$D$. The edge $aw^{\prime\prime}$ joins $G[N[w]]$ to the subgraph induced by the
closed neighbourhoods of the previously chosen vertices. Hence $G[N[D]]$ is
connected. Now let $J$ be a minor of $R$, and contract the connected graph
$G[N[D]]$ into a single vertex. This vertex is adjacent to every branch set of a
minor model of $J$ in $R$, so we obtain a $K_{1}\vee J$ minor of $G$. This
proves (L2).
\end{proof}

\begin{lemma}
\label{Lemma_selector}Every planar graph has a UM-selector of width $3$.
\end{lemma}

\begin{proof}
For $k\in\{1,2,3\}$, let $\mathcal{G}_{k}$ be the class of linear forests
(disjoint unions of paths) if $k=1$, the class of outerplanar graphs if $k=2$,
and the class of planar graphs if $k=3$. We prove by induction on $k$ that
every graph $G\in\mathcal{G}_{k}$ has a UM-selector of width $k$. Since
UM-selectors of components combine, we may assume that $G$ is connected; the
empty graph has the empty selector.

For $k=1$, the graph $G$ is a path $v_{1}v_{2}\cdots v_{n}$. If $n\equiv0,2\pmod
3$, let $S$ consist of the vertices $v_{j}$ with $j\equiv2\pmod 3$. If
$n\equiv1\pmod 3$, let $S$ consist of the vertices $v_{j}$ with $j\equiv1\pmod
3$. In both cases $S$ is independent and every closed neighbourhood contains
exactly one vertex of $S$, so $\varphi\equiv1$ gives a UM-selector of width $1$.

Now let $k\in\{2,3\}$, and take $D$ and $R$ from Lemma \ref{Lemma_JLS}. We first
show that $R\in\mathcal{G}_{k-1}$. Let $J$ be a component of $R$. By (L2),
$K_{1}\vee J$ is a minor of $G$. If $k=3$, then $K_{1}\vee J$ is planar, so $J$
is outerplanar. If $k=2$, then $K_{1}\vee J$ is outerplanar, so it has neither a
$K_{4}$-minor nor a $K_{2,3}$-minor \cite{ChartrandHarary}. Hence $J$ has no
cycle, since otherwise $K_{1}\vee J$ would have a $K_{4}$-minor. Also, $J$ has no
vertex of degree at least $3$, since otherwise $K_{1}\vee J$ would contain
$K_{1}\vee K_{1,3}\supseteq K_{2,3}$. So $J$ is a path, and $R\in\mathcal{G}_{k-1}$
also in this case. By the induction hypothesis, $R$ has a UM-selector
$(S_{R},\varphi_{R})$ of width $k-1$. We set $S=D\cup S_{R}$, $\varphi=\varphi_{R}$
on $S_{R}$, and $\varphi(w)=k$ for every $w\in D$.

The set $S$ is independent, since $D$ is independent by (L1), $S_{R}$ is
independent, and there is no edge between $D$ and $R=G-N[D]$. Now let $x\in
V(G)$. If $x\in N[D]$, then $N[x]$ contains exactly one vertex of $D$, because
two such vertices would be at distance at most $2$, contrary to (L1). So exactly
one vertex of $N[x]\cap S$ has colour $k$, and all other vertices of $N[x]\cap S$
have smaller colours. If $x\in V(R)$, then $N[x]$ contains no vertex of $D$, so
$N_{G}[x]\cap S=N_{R}[x]\cap S_{R}$, and the maximum is unique by the choice of
$(S_{R},\varphi_{R})$. Hence $(S,\varphi)$ is a UM-selector of width $k$ in $G$.
\end{proof}

All our upper bounds place the colours of a UM-selector above the colours of
the other vertices, as follows.

\begin{lemma}
\label{Lemma_combination}Let $(S,\varphi)$ be a UM-selector of width $p$ in a
graph $G$, and let $\psi$ be a colouring of $G-S$ with colours
$\{1,\ldots,q\}$. Define $\sigma=\psi$ on $V(G)\setminus S$ and $\sigma=q+\varphi$
on $S$. Then $\sigma$ is an iUMc-colouring with at most $p+q$ colours. If $G$ has
no isolated vertices and the maximum colour of $\psi$ on $N(s)$ is attained
exactly once for every $s\in S$, then $\sigma$ is also an iUMo-colouring.
Moreover, if $\psi$ is a proper colouring of $G-S$, then $\sigma$ is proper.
\end{lemma}

\begin{proof}
If $\psi$ is proper, then so is $\sigma$, because $S$ is independent and the two
palettes are disjoint. Let $v\in V(G)$. The set $N[v]\cap S$ is non-empty, and
its colours are larger than all colours outside $S$. So the maximum colour on
$N[v]$ is $q+\max\varphi(N[v]\cap S)$, and it is attained exactly once. If
$v\notin S$, then $N(v)\cap S=N[v]\cap S$, so the same holds for $N(v)$. If $v\in
S$, then $N(v)$ lies outside $S$, because $S$ is independent. Thus $\sigma=\psi$
on $N(v)$, and the maximum colour on $N(v)$ is unique by the assumption on
$\psi$.
\end{proof}

\begin{theorem}
\label{Tm_upper}For every planar graph $G$ it holds that $\iUMc(G)\leq4$ and
$\pUMc(G)\leq7$.
\end{theorem}

\begin{proof}
By Lemma \ref{Lemma_selector}, $G$ has a UM-selector $(S,\varphi)$ of width $3$.
With $\psi\equiv1$, Lemma \ref{Lemma_combination} gives an iUMc-colouring with
$3+1=4$ colours. By the Four Colour Theorem \cite{AppelHaken,AppelHakenKoch},
$G-S$ also has a proper colouring $\psi$ with colours $\{1,2,3,4\}$, and with
this $\psi$, Lemma \ref{Lemma_combination} gives a pUMc-colouring with $3+4=7$
colours.
\end{proof}

The bound $\iUMc(\mathcal{P})\geq4$ is proved in \cite{FLRS}. It also follows
from \cite{SedlarSkrekovski}, since every iUMc-colouring is conflict-free on
closed neighbourhoods. Hence we obtain the following.

\begin{corollary}
\label{Cor_iUMc}It holds that $\iUMc(\mathcal{P})=4$, i.e.\ Conjecture \ref{Con_FLRS}(i) is true.
\end{corollary}

\section{\texorpdfstring{A planar graph with $\pUMc=7$}{A planar graph with chi\_pUMc=7}}

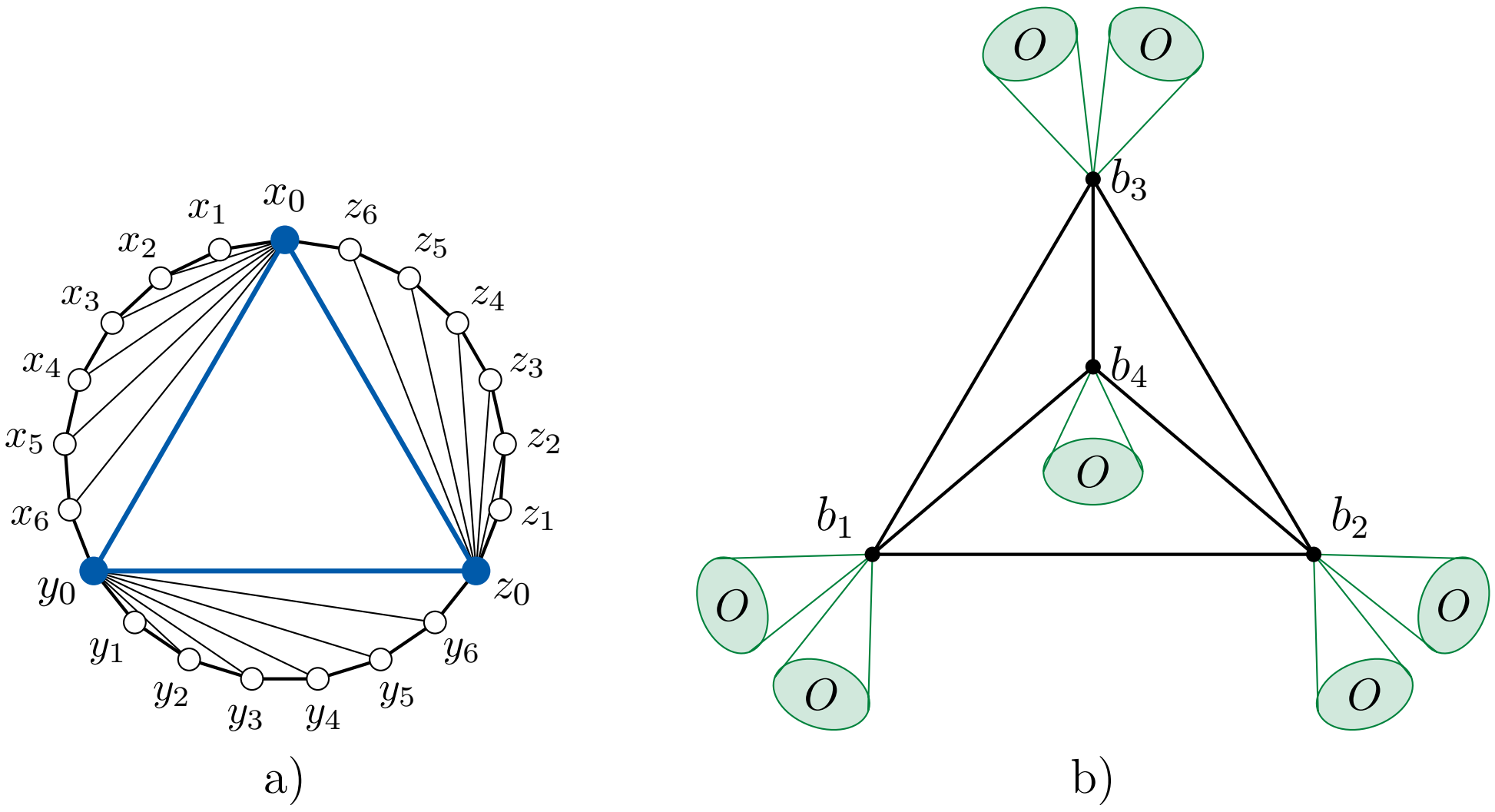
\begin{figure}[ht]
\begin{center}
\begin{tikzpicture}[scale=1.9]
\foreach \i in {0,...,20} {\coordinate (v\i) at ({90+\i*360/21}:1);}
\foreach \i [evaluate=\i as \j using {int(mod(\i+1,21))}] in {0,...,20}
  {\draw[thick] (v\i) -- (v\j);}
\foreach \h/\a in {0/7,7/14,14/0}{
  \foreach \k in {2,...,6}{\pgfmathtruncatemacro{\t}{\h+\k}\draw (v\h) -- (v\t);}
  \draw[very thick,gA] (v\h) -- (v\a);}
\foreach \i in {0,...,20} {\fill[white] (v\i) circle (0.05); \draw (v\i) circle (0.05);}
\foreach \i in {0,7,14} {\fill[gA] (v\i) circle (0.065);}
\foreach \L/\h in {x/0,y/7,z/14}{
  \node[font=\small] at ({90+\h*360/21}:1.19) {$\L_{0}$};
  \foreach \k in {1,...,6}{\pgfmathtruncatemacro{\t}{\h+\k}
    \node[font=\footnotesize] at ({90+\t*360/21}:1.18) {$\L_{\k}$};}}
\node at (0,-1.45) {a)};
\end{tikzpicture}
\qquad
\begin{tikzpicture}[scale=0.95]
\coordinate (b1) at (0,0); \coordinate (b2) at (4,0); \coordinate (b3) at (2,3.4);
\coordinate (b4) at (2,1.7);
\draw[thick] (b1)--(b2)--(b3)--(b1) (b4)--(b1) (b4)--(b2) (b4)--(b3);
\foreach \b/\ang/\r in {b1/200/1.35,b1/250/1.35,b2/-20/1.35,b2/-70/1.35,
  b3/65/1.35,b3/115/1.35,b4/270/0.95}{
  \begin{scope}[shift={(\b)},rotate=\ang]
  \fill[gB!18] (\r,0) ellipse (0.3 and 0.45);
  \draw[gB] (0,0) -- (\r,0.45) (0,0) -- (\r,-0.45) (\r,0) ellipse (0.3 and 0.45);
  \node[font=\small] at (\r,0) {$O$};
  \end{scope}}
\foreach \b/\l/\p in {b1/b_{1}/above left,b2/b_{2}/above right,b3/b_{3}/right,
  b4/b_{4}/right}{\fill (\b) circle (0.07); \node[\p] at (\b) {$\l$};}
\node at (2,-2.05) {b)};
\end{tikzpicture}
\end{center}
\caption{a) The maximal outerplanar graph $O=O_{\mathrm{pUMc}}$ of \cite{FLRS},
drawn with its outer cycle on a circle. The hubs $x_{0}$, $y_{0}$ and $z_{0}$
are the filled vertices and the hub triangle is drawn bold. Every hub is
adjacent to all vertices of the path which starts after it on the circle and
ends with its successor. b) A schematic drawing of the planar graph $G_{7}$ from
Theorem \ref{Tm_lower7}. Every ellipse represents a copy of $O$ whose $21$
vertices are all joined to the corresponding vertex $b_{i}$; the vertices
$b_{1},b_{2},b_{3}$ carry two gadgets each, and $b_{4}$ carries one gadget.}%
\label{Fig_O}%
\end{figure}

Let $O$ be the outerplanar graph $O_{\mathrm{pUMc}}$ of \cite{FLRS}, which is
shown in Figure \ref{Fig_O}(a). It consists of the triangle $x_{0}y_{0}z_{0}$ and
three fans, and we call $x_{0}$, $y_{0}$ and $z_{0}$ \emph{hubs}. The hub $x_{0}$
is adjacent to every vertex of the path $x_{1}x_{2}\cdots x_{6}y_{0}$, the hub
$y_{0}$ to every vertex of the path $y_{1}y_{2}\cdots y_{6}z_{0}$, and the hub
$z_{0}$ to every vertex of the path $z_{1}z_{2}\cdots z_{6}x_{0}$. The graph $O$
is maximal outerplanar and has $21$ vertices and $39$ edges. Its outer cycle is
$x_{0}x_{1}\cdots x_{6}y_{0}y_{1}\cdots y_{6}z_{0}z_{1}\cdots z_{6}x_{0}$. In
\cite{FLRS} it is shown that $\pUMc(O)=5$. For a hub $h$, the next hub in the
cyclic order $(x_{0},y_{0},z_{0})$ is the \emph{successor} of $h$. We denote the
path of the fan of $h$ by $p_{1}p_{2}\cdots p_{7}$, where $p_{7}$ is the
successor of $h$. Note that $p_{2},\ldots,p_{6}$ have degree $3$ in $O$, namely
$N_{O}(p_{i})=\{p_{i-1},p_{i+1},h\}$.

Let $b$ be a vertex of a graph $G$. \emph{Attaching a gadget at} $b$ means adding
a new copy of $O$ to $G$ and joining $b$ to all $21$ vertices of the copy; we
then say that $b$ \emph{carries} this gadget. The vertex $b$ and the copy induce
the planar graph $K_{1}\vee O$, and $b$ is a cut vertex which separates the copy
from the rest of the graph. Hence attaching gadgets preserves planarity.
Moreover, $N_{G}[y]=N_{O}[y]\cup\{b\}$ for every vertex $y$ of the copy. The key
property of the gadget is the following.

\begin{lemma}
\label{Lemma_gadget}Let $b$ be a vertex of a graph $G$ which carries a gadget,
and let $\sigma$ be a proper colouring of $G$ with $\sigma(b)\leq3$ in which the
maximum colour on $N[y]$ is unique for every vertex $y$ of the gadget. Then
some vertex of the gadget has colour at least $6$.
\end{lemma}

\begin{proof}
Suppose, to the contrary, that all vertices of the gadget have colours at most
$5$. We identify the gadget with $O$, so that $N[y]=N_{O}[y]\cup\{b\}$ for $y\in
V(O)$. Let $a=\sigma(b)$. Since $b$ is adjacent to all vertices of $O$, the
vertices of $O$ use only the four colours of $\{1,\ldots,5\}\setminus\{a\}$. We
call a vertex \emph{big} if its colour is $4$ or $5$, and \emph{small} otherwise;
in particular, $b$ is small. Exactly two of the four available colours are
small, namely the colours of $\{1,2,3\}\setminus\{a\}$. The proof is based on the
following claim.

\medskip

\noindent\textbf{Claim A.} \emph{If a hub }$h$\emph{ is small, then the successor
of }$h$\emph{ is also small.}\smallskip

\noindent Let $p_{1}\cdots p_{7}$ be the path of the fan of $h$. For every
$i\in\{1,\ldots,6\}$, the vertices $p_{i}$, $p_{i+1}$ and $h$ form a triangle, so
$p_{i}$ and $p_{i+1}$ have distinct colours, which are both different from $a$
and $\sigma(h)$. Only one small colour remains for them, so $p_{i}$ and $p_{i+1}$
are not both small. Hence at least three of $p_{1},\ldots,p_{7}$ are big. All
these vertices lie in $N[h]$, while $h$ and $b$ are small. Since the maximum
colour on $N[h]$ is unique, it is $5$, as otherwise the colour $4$ would appear
at least three times. So the colour $5$ appears exactly once in $N[h]$, and at
most one vertex of the path has colour $5$. We also need the following property.

\begin{itemize}
\item[$(\ast)$] If $p_{i}$ is small for some $i\in\{2,\ldots,6\}$, then $p_{i-1}$
or $p_{i+1}$ has colour $5$.
\end{itemize}

\noindent Indeed, $N[p_{i}]=\{p_{i-1},p_{i},p_{i+1},h,b\}$, where $p_{i}$, $h$ and
$b$ are small and $p_{i-1}$, $p_{i+1}$ are big. Since the maximum on $N[p_{i}]$
is unique, $p_{i-1}$ and $p_{i+1}$ cannot both have colour $4$.

Let $p_{j}$ be the vertex of the path with colour $5$, and let $j=0$ if there is
no such vertex. By $(\ast)$, every $p_{i}$ with $i\in\{2,\ldots,6\}$ and
$\left\vert i-j\right\vert \geq2$ is big and different from $p_{j}$, so it has
colour $4$. Since adjacent vertices cannot both have colour $4$, we get $j=4$;
otherwise $p_{5},p_{6}$ (if $j\leq3$) or $p_{2},p_{3}$ (if $j\geq5$) would be
adjacent vertices of colour $4$. In particular, $p_{6}$ has colour $4$. The
successor $p_{7}$ of $h$ is adjacent to $p_{6}$ and different from $p_{4}$, so it
has neither colour $4$ nor colour $5$. Thus $p_{7}$ is small, which proves Claim
A.

\medskip

\noindent The three hubs form a triangle, so they have three distinct colours.
Only two of the available colours are big, so at least one hub is small. By
Claim A, applied twice, all three hubs are small. Then the triangle
$x_{0}y_{0}z_{0}$ uses three distinct small colours, although only two small
colours are available. This contradiction proves the lemma.
\end{proof}

Now let $G_{7}$ be the graph obtained from the complete graph $K_{4}$ on the
vertices $b_{1},b_{2},b_{3},b_{4}$ by attaching two gadgets at each of
$b_{1},b_{2},b_{3}$ and one gadget at $b_{4}$; see Figure \ref{Fig_O}(b). The
graph $G_{7}$ is planar, and it has $4+7\cdot21=151$ vertices and
$6+7\cdot60=426$ edges.

\begin{theorem}
\label{Tm_lower7}The planar graph $G_{7}$ has $\pUMc(G_{7})=7$. Consequently,
$\pUMc(\mathcal{P})=7$, and Conjecture \ref{Con_FLRS}(iii) is false.
\end{theorem}

\begin{proof}
We first prove that $\pUMc(G_{7})\geq7$. Suppose, to the contrary, that $\sigma$
is a pUMc-colouring of $G_{7}$ with colours from $\{1,\ldots,6\}$. Let $b$ be a
vertex which carries a gadget and has $\sigma(b)\leq3$. By Lemma
\ref{Lemma_gadget}, every gadget at $b$ contains a vertex of colour $6$. Since $6$
is the largest colour, it appears at most once in $N[b]$. Hence $b$ carries only
one gadget, and no neighbour of $b$ outside of this gadget has colour $6$.

Consequently, each of $b_{1}$, $b_{2}$, $b_{3}$ has colour at least $4$, because
it carries two gadgets. These three vertices form a triangle, so their colours
are $4$, $5$ and $6$. The vertex $b_{4}$ is adjacent to all of them, so
$\sigma(b_{4})\leq3$, and $b_{4}$ has a neighbour of colour $6$ outside of its
gadget. This contradicts the previous paragraph.

Finally, $\pUMc(G_{7})\leq\pUMc(\mathcal{P})\leq7$ by Theorem \ref{Tm_upper}.
\end{proof}

\section{\texorpdfstring{An upper bound for $\pUMo$}{An upper bound for chi\_pUMo}}

For open neighbourhoods, the colouring $\psi$ in Lemma \ref{Lemma_combination}
must also give a unique maximum on $N(s)$ for every $s\in S$, and these sets lie
entirely outside of $S$. For a graph $G$ and an independent set $S\subseteq
V(G)$, let $q(G,S)$ be the least $q$ such that $G-S$ has a proper colouring
$\psi$ with colours $\{1,\ldots,q\}$ in which the maximum colour of $\psi$ on
$N(s)$ is attained exactly once for every non-isolated vertex $s\in S$. By Lemma
\ref{Lemma_combination}, $\pUMo(G)\leq p+q(G,S)$ for every UM-selector
$(S,\varphi)$ of width $p$ in a graph $G$ without isolated vertices.

To bound $q(G,S)$ for planar graphs, we use a theorem of Wendland
\cite{Wendland}. A \emph{facial unique-maximum colouring} of a plane graph is a
proper colouring of its vertices by positive integers such that, for every
face, the maximum colour appears exactly once among the vertices incident with
this face. In the following theorem, as in \cite{Wendland}, a plane graph is a
loopless graph embedded in the plane which may contain parallel edges (see the
end of Section 1 of \cite{Wendland}).

\begin{theorem}
\label{Tm_Wendland}Every plane graph has a facial unique-maximum colouring with
colours $\{1,\ldots,5\}$.
\end{theorem}

The bound $5$ is tight by \cite{LMS}, where the conjecture of Fabrici and
G\"{o}ring \cite{FabriciGoring} that $4$ colours always suffice was disproved.

\begin{lemma}
\label{Lemma_closure}Let $G$ be a simple plane graph and let $S\subseteq V(G)$ be an
independent set. Then there is a loopless plane multigraph $M$ with the vertex
set $V(G)\setminus S$ which contains $G-S$ and has, for every $s\in S$ with
$d(s)\geq2$, an inner face $f_{s}$ with $V(f_{s})=N(s)$.
\end{lemma}

\begin{proof}
Let $U=V(G)\setminus S$. For every $s\in S$ of degree $d\geq2$, let
$v_{1},\ldots,v_{d}$ be the neighbours of $s$ in clockwise order around $s$. We
add the new edges $v_{i}v_{i+1}$ for $i=1,\ldots,d$, with indices modulo $d$;
for $d=2$ these are two parallel edges $v_{1}v_{2}$. The edge $v_{i}v_{i+1}$ is
drawn inside the angle between $sv_{i}$ and $sv_{i+1}$, close to the path
$v_{i}sv_{i+1}$. Finally, we delete all vertices of $S$, and we call the
resulting graph $M$. Since $S$ is independent, all $v_{i}$ belong to $U$, and
the stars of distinct vertices of $S$ have only vertices of $U$ in common. So the
new edges can be drawn without crossings. Since the neighbours
$v_{1},\ldots,v_{d}$ are distinct, no loop arises. Parallel edges may arise, and
we keep them, since otherwise there might be no face whose vertex set is
exactly $N(s)$. Apart from this, $M$ is the facial closure of \cite{FLRS}.

For $s\in S$ with $d(s)\geq2$, the new edges at $s$ form a closed curve, and the
region inside this curve contains only $s$ and its incident edges. After $s$ is
deleted, this bounded region becomes an inner face $f_{s}$ of $M$ which is
bounded by the cycle $v_{1}v_{2}\cdots v_{d}$. Hence $V(f_{s})=N(s)$.
\end{proof}

\begin{theorem}
\label{Tm_pUMo}For every planar graph $G$ without isolated vertices it holds that
$\iUMo(G)\leq\pUMo(G)\leq8$. Consequently, $7\leq\pUMo(\mathcal{P})\leq8$ and
$5\leq\iUMo(\mathcal{P})\leq8$.
\end{theorem}

\begin{proof}
By Lemma \ref{Lemma_selector}, $G$ has a UM-selector $(S,\varphi)$ of width $3$.
Fix a plane embedding of $G$ and let $M$ be as in Lemma \ref{Lemma_closure}. By
Theorem \ref{Tm_Wendland}, $M$ has a facial unique-maximum colouring $\psi$ with
colours $\{1,\ldots,5\}$. Since $G-S\subseteq M$, the colouring $\psi$ is a proper
colouring of $G-S$. For every $s\in S$ with $d(s)\geq2$, the maximum colour of
$\psi$ on $N(s)=V(f_{s})$ is unique, and for $d(s)=1$ this is trivial. Hence
$q(G,S)\leq5$, and Lemma \ref{Lemma_combination} gives $\pUMo(G)\leq3+5=8$. Every
pUMo-colouring is an iUMo-colouring, so also $\iUMo(G)\leq\pUMo(G)$.

It remains to prove the lower bounds. Every pUMo-colouring is also a
pUMc-colouring \cite{FLRS}, so $\pUMo(G_{7})\geq\pUMc(G_{7})=7$ by Theorem
\ref{Tm_lower7}. On the other hand, colour $b_{1},b_{2},b_{3},b_{4}$ by $7,6,5,4$,
respectively, and every gadget by a proper colouring of $O$ with colours
$\{1,2,3\}$. This is a pUMo-colouring of $G_{7}$. Indeed, for a vertex $y$ of a
gadget at $b_{i}$, the unique maximum colour on $N(y)$ is the colour of $b_{i}$.
The unique maximum on $N(b_{1})$ is the colour $6$ of $b_{2}$, and for $i\geq2$
the unique maximum on $N(b_{i})$ is the colour $7$ of $b_{1}$. Hence
$\pUMo(G_{7})=7$. Finally, the lower bound $\iUMo(\mathcal{P})\geq5$ is proved in
\cite{FLRS}.
\end{proof}

\section{Improper open colourings without vertices of degree~2}

In the improper case, the colouring of the closure in the proof of Theorem
\ref{Tm_pUMo} need not be proper. For simple plane graphs, Fabrici and
G\"{o}ring \cite{FabriciGoring} showed that three colours then suffice, and we
extend their argument to plane multigraphs. In this section, a \emph{plane
multigraph} is a finite loopless multigraph embedded in the plane. For a face
$f$, we denote by $V(f)$ the set of vertices incident with $f$. The outer face is
denoted by $\omega$, and all other faces are \emph{inner} faces. We use the
colours black $<$ blue $<$ red.

Fabrici and G\"{o}ring derive their result from \cite[Lemma 7]{FabriciGoring}
(we use the numbering of the journal version), which is stated for all simple
plane graphs. Its proof glues together colourings of two subgraphs $G_{1}$ and
$G_{2}$ (Cases 2, 4.1 and 4.2 there), but this is valid only if no part of
$G_{2}$ lies in an inner face of $G_{1}$. For example, let $G$ consist of a
triangle $xyx_{2}$, a triangle $x_{2}k_{1}k_{2}$ drawn inside it, and a triangle
$x_{2}k_{3}k_{4}$ drawn outside it, and let $c=$ blue. Case 4.1 of
\cite{FabriciGoring} colours the triangle $xyx_{2}$ with $x$ blue and $y,x_{2}$
black, and it colours the graph induced by $x_{2},k_{1},\ldots,k_{4}$ with one of
$k_{1},k_{2}$ blue. Then the face of $G$ incident with $x,y,x_{2},k_{1},k_{2}$
has two blue vertices and no red vertex. We therefore give a complete proof of
the following version of \cite[Lemma 7]{FabriciGoring}, which also covers
multigraphs. It follows the strategy of \cite{FabriciGoring}, but it splits the
graph along closed curves through the outer face.

\begin{lemma}
\label{Lemma_FG}Let $G$ be a plane multigraph in which every component is
incident with the outer face $\omega$. Let $xy$ be an edge incident with
$\omega$, and let $c\in\{\mathrm{black},\mathrm{blue}\}$. Then $G$ has a
colouring with black, blue and red such that \emph{(F1)} $x$ has colour $c$,
\emph{(F2)} $y$ is black, \emph{(F3)} no edge joins two blue vertices,
\emph{(F4)} no vertex incident with $\omega$ is red, and \emph{(F5)} for every
inner face $f$ with $\left\vert V(f)\right\vert \geq3$ the maximum colour on
$V(f)$ is attained exactly once.
\end{lemma}

\begin{proof}
We proceed by induction on $\left\vert V(G)\right\vert +\left\vert
E(G)\right\vert $, and we call a colouring with (F1)--(F5) \emph{good} for
$(G,xy,c)$. We use the following observation, which holds for every plane
multigraph $G$ in which every component is incident with $\omega$.

\begin{itemize}
\item[(R)] Let $z$ be a vertex or an edge of $G$ incident with $\omega$. Deleting
$z$ merges $\omega$ and all faces incident with $z$ into the outer face of $G-z$,
and it does not change the other faces. Hence the vertices incident with the
outer face of $G-z$ are those of $V(\omega)\cup\bigcup_{f\ni z}V(f)$ other than
$z$, and the inner faces of $G-z$ are exactly the inner faces of $G$ not incident
with $z$, with the same vertex sets. Moreover, every component of $G-z$ is
incident with its outer face, since it is a component of $G$ or it contains a
neighbour of $z$ (if $z$ is a vertex) or an end of $z$ (if $z$ is an edge).
\end{itemize}

\smallskip\noindent\emph{Case 1: $G$ is disconnected.} No component of $G$ lies
in an inner face of another component. Therefore every inner face of $G$ is an
inner face of exactly one component, with the same vertex set, and a vertex is
incident with $\omega$ if and only if it is incident with the outer face of its
component. So we may colour the components separately. We colour the component
containing $xy$ by induction for $(xy,c)$, every other component with an edge by
induction for an edge incident with its outer face and the colour black, and
every single-vertex component black.

\smallskip\noindent From now on $G$ is connected, and $W$ denotes the facial walk
of $\omega$.

\smallskip\noindent\emph{Case 2: some vertex $x_{2}$ occurs at least twice in the
closed walk $W$, or some edge $e=x_{2}y_{2}$ not incident with $\omega$ has both
ends in $V(\omega)$.} We split $G$ along a closed curve $\gamma$. In the first
case, $\gamma$ is a simple closed curve which meets $G$ only in $x_{2}$ and
otherwise lies in $\omega$. It leaves $x_{2}$ through one angle of $\omega$ at
$x_{2}$ and returns through a different one, so it splits the edges at $x_{2}$
into two non-empty sets, one in each open region bounded by $\gamma$. In the
second case, $\gamma$ consists of $e$ and an arc from $y_{2}$ to $x_{2}$ through
$\omega$. The two faces at $e$ are inner faces. Since there are no loops, each of
them is incident with an edge other than $e$, and this edge lies in the
corresponding open region bounded by $\gamma$. In both cases, every inner face
of $G$ and the interior of every edge of $G$ not on $\gamma$ are disjoint from
$\gamma$, so each of them lies in one of the two open regions bounded by
$\gamma$. In particular, in the second
case the two faces at $e$ are distinct, since they lie on different sides of
$\gamma$.

Let $G_{1}$ consist of the vertices and edges of $G$ in the closed region bounded
by $\gamma$ which contains $xy$, and let $G_{2}$ consist of those in the other
closed region. Then $G_{1}\cap G_{2}$ is the vertex $x_{2}$ in the first case,
and it consists of $x_{2}$, $y_{2}$ and $e$ in the second case. Both $G_{1}$ and
$G_{2}$ are smaller than $G$. Each $G_{i}$ is connected: since no edge crosses
$\gamma$, a path in $G$ from a vertex of $G_{i}$ to $x_{2}$ stays in $G_{i}$ up to
its first vertex in $G_{1}\cap G_{2}$, and in the second case $y_{2}$ is joined to
$x_{2}$ by $e$. Every inner face of $G$ lies in one of the two open regions, and
its boundary lies in the corresponding $G_{i}$, so it is an inner face of exactly
one $G_{i}$, with the same vertex set. Finally, the vertices incident with the
outer face of $G_{i}$ are exactly the vertices of $G_{i}$ incident with
$\omega$.

We first colour $G_{1}$ by induction for $(xy,c)$. By (F4), the vertices of
$G_{1}\cap G_{2}$ are black or blue. In the second case they are not both blue
by (F3), so after possibly swapping the names of $x_{2}$ and $y_{2}$ we may
assume that $y_{2}$ is black. In the first case, we choose an edge $x_{2}y_{2}$
of $G_{2}$ at $x_{2}$ incident with the outer face of $G_{2}$. Then we colour
$G_{2}$ by induction for $(x_{2}y_{2},c_{2})$, where $c_{2}$ is the colour of
$x_{2}$. The two colourings agree on $G_{1}\cap G_{2}$. Since every inner face of
$G$ and every edge of $G$ other than $e$ belongs to exactly one of $G_{1}$ and
$G_{2}$, their union is good for $(G,xy,c)$.

\smallskip\noindent If Case 2 does not apply, then no vertex occurs twice in $W$.
If $W$ traverses some edge twice, then $W$ traverses only the edge $xy$, so $x$
and $y$ have degree $1$ and $G=K_{2}$. Otherwise $W$ is a cycle $C$ of length at
least $2$ with $xy\in E(C)$. Moreover, every edge with both ends on $C$ belongs
to $C$, again because Case 2 does not apply.

\smallskip\noindent\emph{Case 3: $G=K_{2}$ or $G=C$.} We colour $x$ by $c$ and
all other vertices black, with one exception: if $c$ is black and $G=C$ has
length at least $3$, we colour the neighbour $v\neq x$ of $y$ on $C$ blue. If
$G=C$ has length at least $3$, then $G$ has exactly one blue vertex and no red
vertex. Otherwise $G$ has no inner face with at least three vertices. In both
cases the colouring is good.

\smallskip\noindent\emph{Case 4: otherwise.} Let $e$ be the edge of $C$ at $y$
other than $xy$, and let $\alpha$ be the inner face incident with $e$. If all
vertices of $\alpha$ were on $C$, then $\alpha$ would be bounded by edges of $C$
only, and we would have $G=C$. Hence $\alpha$ is incident with a vertex $u\notin
V(C)$. Let $G^{\prime}=G-e-u$. We apply (R) first to $e$ and then to $u$, which
is incident with the outer face of $G-e$. Thus $G^{\prime}$ satisfies the
assumptions of the lemma, and its inner faces are the inner faces of $G$ not
incident with $u$, because the faces at $e$ are $\alpha$ and $\omega$. Moreover,
the vertices incident with the outer face of $G^{\prime}$ are those of
$(V(C)\cup\bigcup_{f\ni u}V(f))\setminus\{u\}$.

We colour $G^{\prime}$ by induction for $(xy,c)$ and colour $u$ red. Then
(F1)--(F4) hold, because $u\notin V(C)$, every edge at $u$ has a red end, and $e$
has the black end $y$. Let $f$ be an inner face of $G$ incident with $u$. All
other vertices of $f$ lie on the outer face of $G^{\prime}$, so they are not red
by (F4) for $G^{\prime}$. Hence $u$ is the unique maximum on $V(f)$. All other
inner faces of $G$ are inner faces of $G^{\prime}$, so (F5) holds.
\end{proof}

Lemma \ref{Lemma_FG} is also relevant for Theorem \ref{Tm_Wendland}. The proof of
Theorem \ref{Tm_Wendland} in \cite{Wendland} uses \cite[Lemma 7]{FabriciGoring},
whose proof has the gap described above. Adding an edge inside a face incident
with two components changes no set $V(f)$, so that proof may start with a
connected graph. Then, as far as we can see, every graph to which
\cite[Lemma 7]{FabriciGoring} is applied there has all its components incident
with the outer face. For such graphs, Lemma \ref{Lemma_FG} can be used instead,
because in a simple plane graph every inner face has at least three vertices, so
(F5) is equivalent to conditions (5) and (6) of \cite[Lemma 7]{FabriciGoring}.
Hence, as far as we can see, Theorem \ref{Tm_Wendland}, and thus Theorem
\ref{Tm_pUMo}, is not affected.

\begin{theorem}
\label{Tm_iUMo6}Every planar graph $G$ without vertices of degree $2$ and without
isolated vertices satisfies $\iUMo(G)\leq6$. In particular, $\iUMo(G)\leq6$ for
every planar graph of minimum degree at least $3$.
\end{theorem}

\begin{proof}
By Lemma \ref{Lemma_selector}, $G$ has a UM-selector $(S,\varphi)$ of width $3$.
Fix a plane embedding of $G$, and let $M$ be as in Lemma \ref{Lemma_closure}.
Adding an edge inside a face incident with two components changes no set $V(f)$,
so we may assume that $M$ is connected. If $M$ has no edge, then $d(s)=1$ for
every $s\in S$, and we let $\psi\equiv1$. Otherwise we colour $M$ by Lemma
\ref{Lemma_FG} for an edge $xy$ incident with its outer face and
$c=\mathrm{black}$, and we let $\psi$ be this colouring with black, blue and red
replaced by $1$, $2$ and $3$.

Now let $s\in S$. If $d(s)=1$, the maximum colour of $\psi$ on $N(s)$ is
trivially unique. Otherwise $d(s)\geq3$ by assumption, and $f_{s}$ is an inner
face of $M$ with $V(f_{s})=N(s)$ and $\left\vert V(f_{s})\right\vert \geq3$. So
the maximum colour of $\psi$ on $N(s)$ is unique by (F5). Hence Lemma
\ref{Lemma_combination} gives an iUMo-colouring of $G$ with $3+3=6$ colours.
\end{proof}

The assumption on vertices of degree $2$ is needed in this proof. For $s\in S$
with $N(s)=\{v_{1},v_{2}\}$ we need $\psi(v_{1})\neq\psi(v_{2})$, and Lemma
\ref{Lemma_FG} does not guarantee this. For example, let $G$ be obtained from
$K_{4}$ by adding, for every edge $xy$, a new vertex adjacent to $x$ and $y$. The
proof of Lemma \ref{Lemma_selector} may output the six new vertices as $S$, if we
take for $D$ the two new vertices of two disjoint edges of $K_{4}$. Then $\psi$
must be a proper colouring of $K_{4}$, which needs four colours.

\section{Concluding remarks}

In this note we proved that
\[
\iUMc(\mathcal{P})=4,\quad\pUMc(\mathcal{P})=7,\quad7\leq\pUMo(\mathcal{P})\leq
8,\quad5\leq\iUMo(\mathcal{P})\leq8.
\]
This settles parts (i) and (iii) of Conjecture \ref{Con_FLRS} and improves the
bounds of \cite{FLRS}.

We now explain why our method does not give $7$ for $\pUMo$. By Lemma
\ref{Lemma_combination}, it would give $7$ if every planar graph $G$ had a
UM-selector $(S,\varphi)$ of width $3$ with $q(G,S)\leq4$. The proof of the bound $3+4$ for proper
conflict-free colourings in \cite{JLS} does not carry over, since the colour
which it makes unique in $N(s)$ need not be the maximum colour on $N(s)$.
Moreover, $q(G,S)\leq4$ fails for some selectors produced by the construction of
Lemma \ref{Lemma_selector}. To see this, let $G_{0}$ be the plane graph of
Lidick\'y, Messerschmidt and \v{S}krekovski \cite{LMS}, which disproves the
conjecture of Fabrici and G\"{o}ring. This graph has no facial unique-maximum
$4$-colouring, and it has exactly three non-triangular faces; these faces are
pairwise vertex-disjoint and together contain all vertices of $G_{0}$. Let
$G^{\ast}$ be obtained from $G_{0}$ by inserting a new vertex into each of the
three non-triangular faces and joining it to all vertices of that face, and let
$S$ be the set of the three new vertices. The two new vertices in the inner faces
are at distance exactly $3$ from the new vertex in the outer face. So the
construction of Lemma \ref{Lemma_JLS} may start with the latter vertex and then
add the other two. Since $N[S]=V(G^{\ast})$, it stops with $D=S$ and $R$ empty.
With $\varphi\equiv3$, $(S,\varphi)$ is a
UM-selector of width $3$ in $G^{\ast}$. However, $q(G^{\ast},S)=5$. Indeed, a
proper $4$-colouring of $G_{0}$ with unique maxima on the three sets $N(s)$ would
be a facial unique-maximum $4$-colouring of $G_{0}$, which does not exist by
\cite{LMS}, and $q(G^{\ast},S)\leq5$ by the proof of Theorem \ref{Tm_pUMo}. On the
other hand, in all our experiments the first selector tried already satisfied
$q(G,S)\leq4$. We tested all triangulations with at most $12$ vertices, all
connected planar graphs with at most $8$ vertices, and several hundred random
planar graphs with up to $300$ vertices.

Parts (ii) and (iv) of
Conjecture \ref{Con_FLRS} remain open; for part (ii), we do not even know whether
$\iUMo(\mathcal{P})\leq6$. By Theorem \ref{Tm_pUMo}, the bound in part (iv) would
be tight. The smallest order of a planar graph with $\pUMc=7$ is also unknown.
Since all connected planar graphs with at most $10$ vertices have $\pUMc\leq5$
(by computer), this order lies between $11$ and $151$.

\bigskip\noindent\textbf{Computations.} The programs used for the computer
experiments mentioned in this note are available from the authors on request.

\medskip\noindent\textbf{AI declaration.} Anthropic Claude.ai was used in the
preparation of this manuscript. It assisted with the arguments and their
exposition, and with the verification of the proofs. All results were verified
by the authors, who take full responsibility for the content.

\medskip\noindent\textbf{Acknowledgements.} R.~\v{S}krekovski acknowledges
partial support from the Slovenian Research and Innovation Agency program
P1--0383 and the bilateral Slovenian--Croatian project BI-HR/25-27-004.
J.~Sedlar acknowledges support from the
Croatian Ministry of Science, Education and Youth through the bilateral
Croatian--Slovenian project 2025--26, and from the NextGenerationEU instrument
via IP-UNIST-17 (GEORAZ).

\end{document}